\documentclass[a4paper]{article}
\usepackage{amsmath}
\usepackage{amssymb}
\usepackage{amsthm}
\usepackage[english]{babel}
\usepackage[utf8]{inputenc}
\usepackage[T1]{fontenc}
\usepackage{microtype}
\usepackage{enumerate}
\usepackage{mathrsfs}
\usepackage{mathtools}
\usepackage{graphicx}
\usepackage{float}
\usepackage{dsfont}
\usepackage[colorinlistoftodos,prependcaption]{todonotes}
\usepackage{listings}
\usepackage{booktabs}
\usepackage{mathtools}
\usepackage{adjustbox}
\usepackage{wrapfig}
\usepackage{xargs}
\usepackage{xcolor}
\usepackage{csquotes}

\theoremstyle{plain}

\newtheorem*{thm*}{Theorem}

\newtheoremstyle{named}{}{}{\itshape}{}{\bfseries}{.}{.5em}{\thmnote{#3}}
\theoremstyle{named}

\theoremstyle{remark}

\theoremstyle{definition}

\newtheorem*{defi*}{Definition}

\renewcommand{\phi}{\varphi}
\newcommand{\eps}{\varepsilon}

\newcommand{\RR}{\mathbb{R}}

\newcommand\abs[1]{\left| #1 \right|}

\newcommandx{\checkk}[2][1=]{\todo[linecolor=red,backgroundcolor=red!25,bordercolor=red,#1]{check: #2}}
\newcommandx{\fix}[2][1=]{\todo[linecolor=blue,backgroundcolor=blue!25,bordercolor=blue,#1]{fix: #2}}
\newcommandx{\improve}[2][1=]{\todo[linecolor=green,backgroundcolor=green!25,bordercolor=green,#1]{improve: #2}}

\usepackage{import}
\usepackage{libertine}
\usepackage[libertine]{newtxmath}
\usepackage{enumitem}
\usepackage{dsfont}
\usepackage{geometry}

\usepackage{authblk} 

\makeatletter
\DeclareFontFamily{OMX}{MnSymbolE}{}
\DeclareSymbolFont{MnLargeSymbols}{OMX}{MnSymbolE}{m}{n}
\SetSymbolFont{MnLargeSymbols}{bold}{OMX}{MnSymbolE}{b}{n}
\DeclareFontShape{OMX}{MnSymbolE}{m}{n}{
	<-6>  MnSymbolE5
	<6-7>  MnSymbolE6
	<7-8>  MnSymbolE7
	<8-9>  MnSymbolE8
	<9-10> MnSymbolE9
	<10-12> MnSymbolE10
	<12->   MnSymbolE12
}{}
\DeclareFontShape{OMX}{MnSymbolE}{b}{n}{
	<-6>  MnSymbolE-Bold5
	<6-7>  MnSymbolE-Bold6
	<7-8>  MnSymbolE-Bold7
	<8-9>  MnSymbolE-Bold8
	<9-10> MnSymbolE-Bold9
	<10-12> MnSymbolE-Bold10
	<12->   MnSymbolE-Bold12
}{}
\let\llangle\@undefined
\let\rrangle\@undefined
\DeclareMathDelimiter{\llangle}{\mathopen}%
{MnLargeSymbols}{'164}{MnLargeSymbols}{'164}
\DeclareMathDelimiter{\rrangle}{\mathclose}%
{MnLargeSymbols}{'171}{MnLargeSymbols}{'171}
\makeatother

\def\R{\mathbb R}
\def\S{\mathbb S}
\newcommand{\sF}{\mathcal{F}}
\newcommand{\bP}{\mathbb{P}}

\newcommand{\bR}{\mathbb{R}}
\newcommand{\bE}{\mathbb{E}}
\newcommand{\sP}{\mathcal{P}}

\newcommand{\ud}{\mathrm{d}}
\newcommand{\ind}{\mathds{1}}

\theoremstyle{plain}

\date{\today}

\author[1]{Daria Ghilli\thanks{daria.ghilli@unipv.it}}
\author[2]{Cristiano Ricci\thanks{cristiano.ricci@ec.unipi.it}}
\author[3]{Giovanni Zanco\thanks{gzanco@luiss.it}}
\affil[1]{Dipartimento di Scienze Economiche e Aziendali, Università di  Pavia, Via S. Felice Al Monastero 5, Pavia (Italy)}
\affil[2]{Dipartimento di Economia e Management, Università di Pisa, Via Cosimo Ridolfi 10, Pisa (Italy)}
\affil[3]{Dipartimento di Economia e Finanza, LUISS Guido Carli, Viale Romania 32, Roma (Italy)}

\title{A Mean Field Game model for COVID-19 with human capital accumulation}

\begin{document}
\maketitle
\begin{abstract}
In this manuscript we present several possible ways of modeling human capital accumulation during the spread of a disease following an agent based approach, where agents behave maximizing their intertemporal utility. We assume that the interaction between agents is of mean field type, yielding a Mean Field Game description of the problem. We discuss how the analysis of a model including both the mechanism of change of species from one epidemiological state to the other and an optimization problem for each agent leads to an aggregate behavior that is not easy to describe, and that sometimes exhibits structural problems. Therefore we eventually propose and study numerically a SEIRD model in which the rate of infection depends on the distribution of the population, given exogenously as the solution to the the Mean Field Game system arising as the macroscopic description of the discrete multi-agent economic model for the accumulation of human capital. Such model arises in fact as a simplified but tractable version of the initial one.\\

\end{abstract}

\section*{Acknowledgments}
Daria Ghilli, Cristiano Ricci and Giovanni Zanco have been supported by the PRIN project "The Time-Space Evolution of Economic Activities: Mathematical Models and Empirical Applications".\\
Daria Ghilli has been supported by the INdAM-GNAMPA project "Modelli MFGs in Economia per lo studio della dinamica del capitale umano con spillovers spaziali".\\
Giovanni Zanco has been supported by the INdAM-GNAMPA project ""Sistemi con interazione spaziale: convergenza, controllo e applicazioni".

\section{Introduction}
Typical models in epidemiology, as well as in economics, describe the time-evolution of variables - for example, the number (absolute or relative) of susceptible and infectious persons in SIS models. In the last 20 years (at least) there has been an increasing interest for spatial extensions of such models, in which the evolution of variables occurs not only in time but also across space. This is a major feature in the so-called economic geography literature, and in general in any model in which the dynamics is not supposed to be uniform across space. Particularly suitable for spatial extensions are interaction models, namely those in which the dynamics mechanism originates from some sort of interaction between agents and the intensity and effectiveness (or even the form) of such interaction depend on the relative and/or absolute position of the agents. The interest of these spatial models is twofold, as they can bring insight on the effects of the space-distribution of agents on the evolution of variables or, the other way around, on consequences of interaction on the movement of people across space.

The SIR epidemiological model is particularly meaningful under this viewpoint: in its classical formulation \cite{KeMc1927} the dynamics of the variables $S$ and $I$ are governed by the incidence of the disease, that is the product of the probability that a contact with an infectious individual results in a contagion with the average number of susceptibles individuals any infectious individual gets in contact with. Clearly this average number is a quite rough estimate of the actual number of contacts people have, that may well be different in different places, due to differences in people density.

In the epidemics literature there has been some work in the last years \cite{ Guo2020, PaLe2017, TaHa2022} about extending SIS and SIR models with spatial components, before the pandemics began; this type of investigations have been pushed further more recently with the spread of COVID-19  \cite{CoGaMaRo2020,SyWhNi2021, WoLi2020}, since geographical differences in the evolution of the pandemics (both at large and small spatial scales) suggest that epidemics dynamical models should take into account - at least in some phases of the spread of the disease - the distribution of people and of significant economic variables across space. This can be done in several ways, depending on the mechanism and analysis one is interested in and on the pandemic phase one is observing.

The spatial density of people represents a significant factor in the dynamics of diseases; indeed many measures that regulators have taken to oppose the diffusion of the pandemics can be seen as density-dependent (e.g. the requirement to wear masks in crowded places; even lockdown can be seen as a way to drastically reduce large peaks in the spatial density of people).

The attempt to incorporate spatial features into SIR-type models poses two issues: on the epidemiological and economic side, to understand what type of space-dependency might be appropriate; on the mathematical side, passing from an ODE description to a PDE description, to have a tractable model, at least numerically.
A strand of the literature considers space movement of people as time flows by adding diffusion terms, typically a Laplacian or nonlocal diffusion terms, in the equations of SIR-type modes, thus studying a system of PDEs. We refer  e.g. to \cite{FaHo2018,  LiYa2014, Re2004}, and further references in the introductions of \cite{PaLe2017, TaHa2022}.\\
Part of the literature focuses on the way the transition rate depends on the distribution. For example in \cite{PaLe2017} the density is  constant in time but space dependent and arises in the standard true mass incidence function replacing  the total number of individual while in \cite{HuNiEc2013} the density is constant and appears also in the transmission rate with a nonlinear dependence. In both the above mentioned papers the density is exogenous.

In the present paper we are interested in exploiting the dependence on the spatial distribution in the incidence of the disease.
We choose SEIRD  dynamics to model the epidemics, where the population is divided into the compartments: \textit{susceptible, exposed, infectious, removed, dead}. Especially after the COVID19, diverse epidemiological models with additional compartments, such as vaccinated, hospitalized, isolated at home, or age-structured models have been proposed (see e.g. \cite{Te2020} for a mean field game perspective on COVID-19 that considers many of these aspects and \cite{FAGoZa} and its references for age-structured models). Our approach could be adapted to such epidemiological models with no mayor difficulties (but with more consistent modifications in the case of age structured models with age as a variable). We decided to focus on the SEIRD since at the same time it captures the main feature of the COVID19 and it stays simple enough to focus on the main interest of the present paper, that is, the dependence on the distribution.

The main novelty  consists in how we choose  the spatial distribution of agents. The models we propose consider simultaneously the space-time evolution of the  epidemics and of a key economic variable, the human capital. 
The description of the evolution of human capital is of crucial importance in economics and has attracted more and more interest nowadays, see \cite{BoDiAz2008} for a review of human capital accumulation models on an epidemic setting, or \cite{Blea2010} where some empirical evidence is proposed. In particular, in the model we consider for the space-time evolution of human capital, the dynamics of an individual's human capital are affected by the human capital of nearby individuals through a spatial interaction term, that is, we account for \textit{spatial spillover} on accumulation of human capital. Individuals control their position in space to maximize their gain from such spillover, considering at the same time that moving towards areas with high human capital is a costly action.  Assuming that no individual has overwhelming influence on the system with respect to the others, optimization of each individual's position gives rise to a symmetric game among economic agents (the individuals).
 When the number of the agents is large 
 the aggregate dynamics are typically described by a Mean Field Game (MFG hereinafter) system, arising heuristically as the limit as the number of agents tends to infinity. The typical MFG system consists of  two partial differential equations: an Hamilton-Jacobi-Bellman (HJB hereinafter) equation describing the optimal control problem of the agents, and a Fokker-Planck (FP hereinafter) equation, describing the  time-distribution of the population with respect to space and human capital.  
An important motivation for the MFG approach we choose for our economy lies in the fact that a solution to the MFG system is expected to provide a good approximation of Nash equilibria for the discrete multi-agent model when the number of agents is large. Moreover, it allows for substantial computational simplifications. In fact, it is well known that solving $N$ coupled optimization problems, or analogously an optimization problem of dimension $N$, becomes unviable from the numerical point of view as soon as $N$ gets moderately large ($N = 10$ is already beyond reach in the general case). On the contrary the problem in the mean field limit 
reduces to a system of two coupled partial differential equations that can be handled in a simpler way.

For a detailed presentation of such model - without reference to epidemics - and in particular of the mathematical framework, we refer to the forthcoming paper \cite{GRZ22}. 
The MFG model is proposed and studied analytically and numerically in \cite{GRZ22} and the details of the model are discussed below in Section \ref{sec:MFG}.

In the present paper we link the MFG dynamics just described for the human capital with epidemics dynamics; we propose several epidemiological models with different aims and mathematical difficulties.  We begin by assuming that the space-time distribution of individuals is endogenously determined by their optimal choices. First we introduce a new micro-economic model for human capital-based interaction of individuals during the spread of a disease that is in our opinion particularly realistic;
in this model individuals tend to increase their level of human capital by being close to other individuals and at the same time to maintain social distancing, in order to minimize the probability of getting infected. 
Our goal is to describe the interplay between the effects on the human capital of each individual caused by the spatial proximity of other individuals and the need of social distancing due to the current level of the epidemics. The corresponding optimization problem is however too difficult to solve, and it is not clear how one should formulate the problem at the MFG level. The major technical difficulties arise since we assume that each class has different goals with respect to keeping social distance and that the epidemics dynamics for individuals is described through piecewise deterministic stochastic processes that jump from a compartment to the next at random times whose intensity depends on the interaction. To the best of our knowledge, the only attempts in this direction available in the epidemic modeling literature are \cite{Doncel2020,Olmez2021,Petrakova2021}, where, however, the optimization problem is the same for all classes, thus avoiding a major technical difficulty of the model. We also mention \cite{Be20}, where MFGs with change of states are studied. 
We underline that these issues are related with the optimization problem only: the appearance of random jumps is not a problem by itself if one wants to understand the aggregate dynamics in the setting of mean field limits. These have been broadly studied, especially in the case of biological applications, see e.g. \cite{catellier2021mean,oelschlager1989derivation}, when however no optimization is involved.
Nonetheless a proper understanding of how to derive the aggregate MFG dynamic for jump processes in the controlled case we propose is still beyond reach for the current knowledge on MFG, even tough the theory of optimal control of jump process has been already developed to many extents, see e.g. \cite{Calvia18,Calvia20} and references therein. We also remark that interacting systems with spatial interaction but without control are also mathematically well understood, see for example \cite{FlPrZa, Zanco} and references therein.
Therefore in the subsequent part of the paper we bring some successive simplifications that allow us to study the aggregate dynamics in the MFG system. The first possibility we consider is to neglect the transition between epidemiological states. In this way we obtain a MFG systems for each compartment $S,I,E,R$ of our epidemiological model 
MFGs with different compartments (often called species in similar models), but where no transitions between compartments are allowed have been proposed for example in  \cite{Achdou2017, Bensoussan2018, Cirant2015, Feleqi2013}. However, we believe that ruling out transition between compartments is too restrictive as assumption. Therefore we formulate a MFG model where we incorporate transitions but we assume that there is no distinction among compartments in the optimization goals. The macroscopic evolution in this case is described by one MFG system and one differential equations for each compartment. Moreover we suppose that the  infection rate depends on the space-time distribution of the agents. In this model it is easily seen that the space-time distribution must satisfy a structural condition, namely to be independent of time. For this reason  we choose for our model the space marginal of the stationary solution to the FP equation in the MFG system. This actually means that to have a model that is sufficiently simple to be tractable, at least numerically, the distribution has to be exogenous.
From the modeling point of view, considering a stationary distribution can be seen as corresponding to the assumption that the mechanism that links the evolution of human capital with the individuals' positions has run for a long time and has stabilized. Incorporating such stationary solution exogenously in our SEIRD model corresponds then to assuming that the epidemics dynamics has negligible effect of the distribution of people; this may be seen as the two distinct dynamics, the one described by the MFG system and the one described by the SEIRD model, happening on different time-scales, with a much shorter time-scale for the epidemics dynamics.
Considering  the space marginal of the solution to the FP equation means that we do not take into account the influence of the accumulation of human capital on the evolution of the epidemics. Up to our knowledge, there is little amount of study or evidence in the literature  of the influence of such variable on the incidence term in epidemiological models, see \cite{BeEbGeJoOi2021, GoLi2020}. On the contrary, some authors have been focusing on the effect of the evolution of epidemics on the level of human capital \cite{Blea2010, DeShLiLi2021}; this phenomenon is considered in the first micro-economic model we propose, but is inevitably ruled out by the simplifications we have to perform.
At the end of these simplifying steps we obtain a model that, omitting some details that are discussed thoroughly in the following sections, takes in the epidemics dynamics the form
\begin{equation*}
\begin{cases}
\partial_{t}S(t,x) &= - \beta\big(\mu(x)\big) \displaystyle\frac{S(t,x)}{\mu(x)}\int K_\chi(x-y) I(y)\ud y \ ,\\[3mm]
\partial_{t}E(t,x) &= \beta\big(\mu(x)\big)\displaystyle \frac{S(t,x)}{\mu(x)}\int K_\chi(x- y) I(y)\ud y - \theta E(t,x)\ ,\\[3mm]
\partial_{t}I(t,x) &= \theta E(t,x) - \lambda  I(t,x) - \delta I(t,x)\ ,\\[2mm]
\partial_{t}R(t,x) &= \lambda I(t,x)\ ,\\[2mm]
\partial_{t}D(t,x) &= \delta I(t,x)\ .
\end{cases}
\end{equation*}
With respect to the classical SIR dynamic, there are the  two additional compartments $E$  (\textit{exposed} individuals) and $D$ (\textit{dead} individuals). The presence of the $E$ compartment implies that there is a latency period, typical of many infectious diseases, during which individuals have been infected but are not yet infectious.  We focus on modeling just deaths due to the infection and we do not consider neither births nor death due to other causes than the infection. Moreover we do not assume diffusion of individuals, as it seems unfit to this type of models at large spatial scales, see also the discussion in the introduction of \cite{PaLe2017}.
On the other hand, we assume that  the infection propagates in space, that is, a susceptible individual can be infected by an infectious individual provided they are close enough. This propagation effect is modeled through a convolution in the space space variable, where the support of the convolution kernel quantifies the threshold for possible contagion to occur. For analogous approaches we refer to \cite{CoGaMaRo2020, TaHa2022} and references therein. We remark that differently from \cite{CoGaMaRo2020} we do not consider age dependence, focusing only on space dependence. As mentioned above, the way the transition rate depends on the distribution is inspired by \cite{HuNiEc2013}.

It is our belief that the models we discuss herein can be also effectively applied to the description of other diseases as well as other economic variables, modifying parameters and related properties of the models accordingly.

In Section \ref{sec:listOfBadModels} we present the first three different models that we just discussed, 
highlighting for each one the mathematical difficulties that are encountered and why these difficulties make those model unviable. 
In Section \ref{sec:model} we present the joint human capital-epidemiological model used in the rest of the paper, that results from the simplifications we explained above. In particular in Subsection \ref{sec:epimodel} the epidemiological component is presented, while in Subsection \ref{sec:MFG} we present the MFG component of the model related to human capital accumulation. Finally in Section \ref{sec:numericaSEIRD} numerical results  for the model given in Section \ref{sec:model} are presented.



\section{Mean Field Game epidemiological models for economic variables}\label{sec:listOfBadModels}
As stated in the introduction we aim to propose Mean Field Game (MFG hereafter) type model which consider simultaneously both the time evolution of economic variables, such as human capital, as well as the time evolution of epidemics. The model we choose to describe the economy (see Section \ref{sec:MFG} for a detailed description) is a spatial model of accumulation of the human capital at agent level, in presence of an epidemic pathogen. Agents are characterized an epidemiological state which represent the stage of infection he is currently on, e.g. \emph{Susceptible}, \emph{Exposed}, \emph{Infected} etc., as well as by a position in space $x$ and by a level of human capital $h$.

We assume that each agent behave optimally in the sense of intertemporal utility maximization, that is each agent is assumed to be \emph{forward looking}. Each agent utility depends not only on the characteristic of the agent itself but also on the configuration of system as a whole, i.e. depends on the position, endowment of human capital and epidemic state, of all other agents. 
The interplay between economic and epidemic variables is expressed in the objective function. Agents tend to maximize their level of human capital by being close (in the space variable) to other agents with a higher level of $h$. This constitutes a form of \emph{positive spatial spillover} with regards to human capital. At the same time, depending on their current epidemiological state, they tend to maintain social distancing, in order to minimize the probability of becoming infected, therefore increasing their probability to die. These two mechanisms represents a form of \emph{competing effect} between the tendency of agglomerating, therefore being close to each other, and the necessity of social distancing. 

In this section we present an overview on different possible mathematical approaches to the problem described above, where different levels of mathematical detail are employed. However, we will see how none of the approaches employed in this section are viable, either because of technical difficulties coming from the mathematical description, or due to lack of some of the key mechanism needed for a proper understanding of the dynamics, which makes these models unsatisfying for our scopes. In Section \ref{sec:model} we carefully detail the actual model that we consider, that derives from some thoughtful consideration of the issues, described in this section. 


\subsection{Change of species and different objective functionals}\label{subsec:dream}
We suppose there are $N$ agents at time $t=0$ in the system; agents are indexed by $j\in\{1,\dots,N\}$. Each agent $j$ is characterized by her epidemiological state (susceptible, exposed, infectious, recovered, dead), her position in a physical region $O$ and her human capital; therefore to each agent $j$ we associate a tuple $(b_j,x_j,h_j)\in\left\{S,E,I,R,D\right\}\times O\times \bR_+$. Agents behave differently depending on their epidemiological state; their goal is to maximize the gain from human capital trying at the same time not to be infected. In order to maximize his level of  capital, each agent chooses at each time $t$ his velocity $v_j(t)$ (intended as speed and direction) the drives his motion across space. Moreover at each time he also maximize his utility depending on the endowment of human capital, by choosing optimally the control $f_j(t)$, that represent the fraction of capital that is reinvested to improve further production. The remaining part $(1-f_j(t))$ appears directly into instantaneous utility, and represent the fraction of human capital that is not reinvested and is instead enjoyed instantaneously at each time $t$.

The $(x,h)$ dynamics of the model is
\begin{equation}
\label{eq:dyn_dream}
\begin{cases}
\ud x_j(t)=\left[v_j(t)\ud t+\epsilon\ud W(t)\right]\ind_{\{S,E,I,R\}}(b_j(t))\\
\ud h_j(t)=\left[f_j(t)h_j(t)^\alpha\bar h_{-j}(t)^\xi\ud t-\zeta h_j(t)\ud t+\sigma h_j(t)\ud W_j(t)\right]\ind_{\{S,E,I,R\}}(b_j(t))\\
x_j(0)=x_{j,0},\quad h_j{0}=h_{j,0},\quad b_j(0)=b_{j,0};
\end{cases}
\end{equation}
clearly $b_j=D$ if and only if the agent is dead because of the infection thus there are no dynamics in this case, hence the indicator functions. The epidemiological state affects the functional to be maximized. In particular each agent in one of the classes $\{S,E,I,R\}$ tries to maximize his level of human capital by either reinvesting his current endowment, or by being close to other individuals with an endowment higher than his own. This mechanism is present in all the epidemiological classes, apart from those in the class $D$ which for whom the utility function is $0$ at each time $t$.

On the contrary, the tendency to keep the distance from other individuals (social distancing) heavily depends on the class each agent belongs to. In particular we assume that susceptibles $(S)$ aim to keep their distance from any other individual in the classes $\{S,E,I,R\}$. The rationale behind this choice is the following: in an epidemic scenario is not always possible to tell if an individual is healthy $(S)$ or has instead been exposed $(E)$. Therefore agents in $(S)$, in absence of perfect knowledge about the state of other individuals will tend to keep social distancing from all other agents, of whom they don't know the epidemiological state. At the same time, in the class of $SEIRD$ models infected individuals $(I)$ are assumed to be symptomatic, therefore they are easily identifiable and are kept at distance as well. Moreover since we don't assume that recovered $(R)$ can get reinfected, and since individuals who are recovered are aware of their status, the mechanism of social distancing of susceptibles does not involve recovered individuals. No interaction is also considered with the dead $(D)$ class since the only disadvantage of belonging to the class $D$ is that utility accumulation stops. 
In particular for the $(S)$ class the tendency to keep separated from other, non infected, individuals is influenced by the value $C_I(t)$. For any given time $t$ $C_I(t)$ represent the percentage of the total population which is infected, which is a function of the true percentage of infected people; thus
\begin{equation*}
 C_I(t)=g\left(\frac{1}{N}\sum_{j=1}^N\ind_{\{I\}}(b_j(t))\right).
\end{equation*}
This expresses the fact that in absence of perfect information on individuals who have been  infected, the only information that is available to public is the percentage of the total population that has been identified (via television, journals of other means of information available to the public). Therefore, the value of $\frac{1}{N}\sum_{j=1}^N\ind_{\{I\}}(b_j(t)))$ is assumed to be known by the population. The effect that this information has on each individual is described by the function $g$. In particular we assume that  $g$ is monotone increasing, $g(0) = 0$ and $g(1) = 1$. Note in particular that the argument of $g$ is bounded between $0$ (no infected) and $1$ (everyone infected) so that the value of $g(1)$ represent the maximum coefficient describing the tendency of agents to avoid each other.

Individuals in the infected class $(I)$, while maintaining the desire to improve their level of human capital in the objective function, tend to keep their distance from any other individual. Notice that the behavior of individuals in $(I)$ is not affected by the coefficient $C_I(t)$. In fact, infected individuals are aware of their own state (since thy are symptomatic), therefore they try to keep their distance from individuals in other classes to avoid to infect them, and not because they are unaware of their epidemiological state. 
Finally, individuals in the recovered class $(R)$ will only keep the distance from infected individuals. Whilst it is true in the model that they cannot get reinfected it does make sense from a real world perspective to avoid contact with infected people.

\begin{multline}
  \label{eq:J_dream}
  J\left(b_{j,0},x_{j,0},h_{j,0};v_j(\cdot),f_j(\cdot)\right)=J^{S,E}\left(x_{j,0},h_{j,0};v_j(\cdot),f_j(\cdot)\right)\ind_{\{S,E\}}(b_j)\\+J^I\left(x_{j,0},h_{j,0};v_j(\cdot),f_j(\cdot)\right)\ind_{\{I\}}(b_j)+J^R\left(x_{j,0},h_{j,0};v_j(\cdot),f_j(\cdot)\right)\ind_{\{R\}}(b_j)
  \end{multline}
where
\begin{multline*}
J^{S,E}\left(x_{j,0},h_{j,0};v_j(\cdot),f_j(\cdot)\right)=\\
\bE\left[\int_0^{+\infty}e^{-\rho t}\Bigg(u\left(\left[(1-f_j(t))h_j(t)^\alpha\right]^{1-\gamma}\bar h_{-j}(t)^\gamma A\left(x_j(t)\right)\right)\Bigg)\ud t\right]\\
-\bE\Bigg[\int_0^{+\infty}e^{-\rho t}\Bigg(a\left(v_j(t)\right)+C_I(t) V_\omega\Bigg(\frac{1}{N}\sum_{k\neq j,k=1}^N\tilde K_\chi\left(x_j(t),x_k(t)\right)\ind_{\{S,E,I,R\}}\left(b_k(t)\right)\Bigg)\ud t\Bigg],
\end{multline*}
\begin{multline*}
J^{I}\left(x_{j,0},h_{j,0};v_j(\cdot),f_j(\cdot)\right)=\\
\bE\left[\int_0^{+\infty}e^{-\rho t}\left(u\left(\left[(1-f_j(t))h_j(t)^\alpha\right]^{1-\gamma} \bar h_{-j}(t)^\gamma A\left(x_j(t)\right)\right)\right)\ud t\right]\\
-\bE\left[\int_0^{+\infty}e^{-\rho t}\left(a\left(v_j(t)\right)+ V_\omega\left(\frac{1}{N}\sum_{k\neq j,k=1}^N\tilde K_\chi\left(x_j(t),x_k(t)\right)\ind_{\{S,E,I,R\}}\left(b_k(t)\right)\right)\right)\ud t\right]\ ,
\end{multline*}
\begin{multline*}
J^{R}\left(x_{j,0},h_{j,0};v_j(\cdot),f_j(\cdot)\right)=\\
\bE\left[\int_0^{+\infty}e^{-\rho t}\left(u\left(\left[(1-f_j(t))h_j(t)^\alpha\right]^{1-\gamma}\bar h_{-j}(t)^\gamma A\left(x_j(t)\right)\right)\right)\ud t\right]\\
-\bE\left[\int_0^{+\infty}e^{-\rho t}\left(a\left(v_j(t)\right)+C_I(t)V_\omega\left(\frac{1}{N}\sum_{k\neq j,k=1}^N\tilde K_\chi\left(x_j(t),x_k(t)\right)\ind_{\{I\}}\left(b_k(t)\right)\right)\right)\ud t\right]\ .
\end{multline*}
In the above equations $u$ is an utility function of consumption ($[(1-f_j(t))h_j(t)^\alpha]^{1-\gamma}$), spatial spillovers on consumption ($\bar h_{-j}(t)^\gamma$) and local amenities ($A(x_j(t))$); $-a(v_j(t))$ represents the cost for the energy employed in the displacement. The function $-V_\omega$ represents a penalization due to being close to other people, thus increasing the probability of getting sick. We assume that $V_\omega(y)=y^\omega/\omega$ with $\omega>1$. The function $\tilde K_\chi(y,z)$ depends only on the distance $d(y,z)$, is positive, decreasing and supported in $[0,\chi]$.

Finally we need to describe the dynamics in the variable $b$; they are pure jump random dynamics that determine the transition $S \to E \to I \to \{R,D\}$. If agent $j$ is susceptible, that is, $b_j(t)=S$, she becomes exposed over an infinitesimal time interval of length $dt$ with probability proportional to 
\begin{equation}
\label{eq:JumpRateStoE}
\phi\cdot \left(\frac{1}{N}\sum_{k\neq j,k=1}^N\tilde K_\chi\left(x_j(t),x_k(t)\right)\ind_{\{I\}}\left(b_k(t)\right)\right) \cdot dt, \quad \phi > 0.
\end{equation}
The rationale behind the previous formula is the following: at each time individual $j$ in the $S$ class has a certain amount of infected people around him. A quantification of the vague expression \emph{around} is given by the function $\tilde{K}_{\chi}$ that measure how is it likely to get infected by staying close to another infected individual. Expression \eqref{eq:JumpRateStoE} is a sum that ranges over all other infected agents in the system, therefore the we are assuming that each individual contributes to the probability of agent $j$ to get infected by
\begin{equation*}
\frac{\phi}{N} \cdot \tilde K_\chi\left(x_j(t),x_k(t)\right)\ind_{\{I\}}\left(b_k(t)\right) \cdot dt
\end{equation*}
over an infinitesimal time interval of length $dt$.

More rigorously, the random clock at which an individual switch from the class $S$ to the class $E$ can be described as the first jump time of inhomogeneous Poisson Process, in which the rate of jump depend itself on the spatial distribution of infected individual (those such that $b_{j}(t) = I$). In particular if $N^{1,j}_{t}$ is standard Poisson Process with rate $\lambda = 1$, and $\tau_{j,1}$ is his first jump time, that is an exponential random variable of mean one, then the time of transition for agent $j$ from the class $S$ to $E$ is defined as 
\begin{equation*}
\tau_{j,S\to E} = \inf\left\{t > 0 \,\Bigg\vert\, \phi \cdot \left(\frac{1}{N}\sum_{k\neq j,k=1}^N\tilde K_\chi\left(x_j(t),x_k(t)\right)\ind_{\{I\}}\left(b_k(t)\right)\right) > \tau^{j,1} \right\}.
\end{equation*}
This corresponds to the first jump time of an inhomogeneous Poisson Process of rate given by \eqref{eq:JumpRateStoE}.

The transition from the class $E$ to $I$ is easier to describe, since the transition rate does not depend on the distribution of all other agents but is assumed constant and exogenous. This is consistent with the fact that the expected time between when an individual is exposed to the disease and when he becomes infectious depends only on the type of the disease. 
If at time $\tau_{j,S\to E}$ individual $j$ get exposed, then he becomes infectious at a random time $\tau_{j,S\to E} + \tau_{j,E\to I}$ where $\tau_{j,E\to I}$ is an exponential random variable with intensity $\tilde{\theta}$. 
If $t_{j,I}=\inf\left\{t\colon b_j(t)=I\right\} = \tau_{j,S\to E} + \tau_{j,E\to I}$ agent $j$ either recovers or dies, depending on which of the random times $t_{j,I}+\tau_{j,I\to R}$ and $t_{j,I}+\tau_{j,I\to D}$ occurs first; here $\tau_{j,I\to R}$ and $\tau_{j,I\to D}$ are exponential random variables with intensity $\tilde{\gamma}$ and $\tilde{\delta}$, respectively; $1/\tilde{\gamma}$ represents the expected waiting time of recovery and $1/\tilde{\delta}$ the expected waiting time of death, thus we assume that $\tilde{\delta}>\tilde{\gamma}$.\\

Our model could be described in terms of piecewise deterministic processes, as discussed for example in \cite{Davis}, but such description would require the introduction of several technical mathematical details that are unnecessary for our discussion in the present paper. As discussed in the introduction, the structure of this model, in particular the presence of the interaction term inside the intensity of the random jumps from compartment $S$ to compartment $I$ and the presence of distinct functionals to be optimized for each compartment, makes it extremely challenging to even formulate a MFG system that may approximate the behavior of individuals as described by this model when their number is large.

\subsection{Mean Field Game with multiple species but no transition between species}\label{subsec:dream0}
In this section we present a simplified version of the previous model. Agents are still characterized by a position in space $x$ and by their endowment of human capital $h$, as well as an epidemiological state $b \in \{S,E,I,R,D\}$. However, as discussed in the previous section, the transition of agents from one class to the other is one of main difficulties that appears in these class of models. Therefore, in the current section, we consider yet a MFG type model, but without allowing agents to transfer from one class to another. 
In this section we don't formulate the problem directly at the agent level, but instead, we focus on the aggregate dynamic. The reasons are twofold: first the microscopic description in the case considered here would be the same of that argued in the previous section, but without considering the random jump times describing transition between classes, i.e. the variable $b \in \{S,E,I,R,D\}$ of each agent would be independent of $t$. Second, since neglecting the change of species actually allow us to understand (at least heuristically) which one is the proper aggregate dynamic, it is more interesting to focus on that directly.  

In the same flavor as in subsection \ref{subsec:dream}, the optimal control problem of each individual is influenced by the state of the pandemic by the penalization term 
\begin{equation}\label{eq:convKchi}
V_\omega(K_\chi*_x \nu^k), \quad k \in \{S,E,I,R\},
\end{equation}
due to being close to the classes $S,E,I$, the classes $S,E,I,R$ and $I$ for susceptible, infected and removed, respectively; here we use the notation
\begin{equation}\label{eq:convdef}
\left(K_\chi*_x\nu^k \right)(t,x,h)=\int_{O}K_\chi(x,y)\nu^k(t,y,h)dy, \quad k \in \{S,E,I,R\}
\end{equation}
where  $K_\chi$ is a positive and decreasing function of the distance between $x, y$ supported in $[0,\chi]$, 
$$
\nu^{S,E}(t,x,h)=C_I(t)\left(\mu^S(t, x, h)+ \mu^E(t,x,h)+\mu^I(t,x,h)+\mu^R(t,x,h)\right),
$$
$$
\nu^I(t,x,h)=\mu^S(t, x, h)+ \mu^E(t,x,h)+\mu^I(t,x,h)+\mu^R(t,x,h), 
$$
$$
\nu^R(t,x,h)= C_I(t)\mu^I(t,x,h)
$$
and $\mu^k$ is the distribution of the  individuals in the compartment $k \in \{S,E,I,R\}$.
The penalization term depends again on the percentage of people infected through the term $C_I$ that now takes the form
\begin{equation}\label{eq:CIdef}
C_{I}(t) := g\left(\int_{O\times \R_+}\mu^I(t,x,h)\,dxdh\right).
\end{equation}

The macroscopic description is given by four MFG systems with solutions couples $(V^k,\mu^k)$, $k \in \{S,E,I,R\}$,  the value function of the class $k$ and the distribution of the population in the class $k$, respectively and a differential equation for $\mu^D$. Each distribution for $k \in \{S, E, I, R\}$ evolves according to a FKP equation describing the optimal evolution, that is whose drift is the optimal control for each class. The value function $V^k$ of each optimal control problem is a solution to an elliptic HJB equation where the HJB equations differ from each other by the penalization terms \eqref{eq:convKchi}.
For the motivation to our choice of the epidemiological model we refer to subsection \ref{sec:epimodel} and for a description of the economic model to subsection \ref{sec:MFG}.\\
More precisely, the MFG systems are the following for each class $k \in \{S,E,I,R\}$: 
$$
\begin{cases}
\rho V^k=H_1\big(x,h,\mu^k(t), \partial_hV^k\big)+ \frac12 \sigma^2h^2 \partial^2_{hh}V^k+\frac12 \eps^2 \partial^2_{xx}V^k+H_0(\partial_xV^k) - V_\omega\left(K_\chi*_x\nu^k\right)\\
\partial_{t}\mu^k= \frac{1}{2}\sigma^2\partial^2_{hh}(h^2\mu^k)+\frac{1}{2}\eps^2\partial^2_{xx}\mu^k-\partial_x\left(\partial_p H_0(\partial_x V^k)\mu^k\right)-\partial_h\left(\partial_p  H_1(x,h,\mu^k(t), \partial_hV^k)\right)\\
\partial_t\mu^D=\delta \mu^I,
\end{cases}
$$
where $H_1$ and $H_0$ are Hamiltonians defined in subsection \ref{sec:MFG}.

The model described above has the advantage of being more manageable from a mathematical point of view with respect to the one proposed in the paragraph \ref{subsec:dream}, but has the remarkable drawback of  not taking into account the change of species. Namely, even if in model discussed in this section each agent acts trying to maximize his own objective functional, depending on the epidemiological class he belongs to, transition between classes are not allowed. For some applications having different classes of agents interacting with each other and each one maximizing his own personal utility depending on the class may be appropriate. However in the framework of an ongoing epidemic the mechanism of transitioning from one class to another is fundamental, therefore cannot be neglected. 
In the following paragraph we propose a more consistent, to our point of view, epidemiological model   where change of species are taken into account but each individual solve the same optimal control problem.

\subsection{Mean Field Game for agent distribution coupled with SEIRD dynamic}\label{subsec:dream00}
The macroscopic description of the model of subsection \ref{subsec:dream0} is relatively simple, at least with respect to the one of subsection \ref{subsec:dream}. However,  neglecting  the transition from one class to another is a strong assumption on the modeling point of view. The goal of this subsection is to present a further simplification of the model described in \ref{subsec:dream} which considers the transition from one class to another. In order to overcome the technical difficulties due to the presence of random jump times describing the transition within classes considered in \ref{subsec:dream0}, we assume the penalization of being close to other people is the same for every class. As in subsection \ref{subsec:dream0}, each agent is characterized by its position $x$, endowment of human capital $h$ and epidemiological state $b \in \{S,E,I,R\}$.  Again and for analogous reasons as in subsection \ref{subsec:dream0}, we consider just the aggregate dynamic. In the limit as number of agents tend to infinity,   the penalization term is
$$
-C_I(t)V_\omega(K_\chi*_x\mu)
$$
where $K_\chi*_x\mu$ is defined in \eqref{eq:convdef} and $C_I(t)$ in \eqref{eq:CIdef}
and
\begin{equation}\label{eq:sumM}
S(t,x,h)+E(t,x,h)+I(t,x,h)+R(t,x,h)+D(t,x,h)=\mu(t,x,h).
\end{equation}
The macroscopic description, as the number of agents tend to infinity, is a  MFG coupled with a SEIRD dynamic. The solution of the MFG is a couple $(V, \mu)$, where  $\mu$ represents the total distribution of the individuals and solves a Fokker-Planck equation with  drift  the optimal control and $V$ is the value function of the optimal control problem and solves an HJB equation.  The evolution of the pandemics is given by five equations one for each class $S,E,I,R,D$. 

The resulting MFG system coupled with the equations describing the evolution of the pandemics is the following system
\begin{equation}\label{MFGnonfunziona}
\begin{cases}
\rho V=H_1\big(x,h,\mu(t), \partial_hV\big)+ \frac12 \chi^2h^2 \partial^2_{hh}V+\frac12 \eps^2 \partial^2_{xx}V+H_0(\partial_xV) - C_{I}(t)V_{\omega}\left(K_\chi*_x\mu\right)\ ,\\[3mm]
\partial_t\mu=\frac{1}{2}\chi^2\partial^2_{hh}(h^2\mu)+\frac{1}{2}\eps^2\partial^2_{xx}\mu-\partial_x\left(\partial_p H_0(\partial_x V)\mu\right)-\partial_h\left(\partial_p  H_1(x,h,\mu(t), \partial_hV)\mu\right)\ ,\\
\partial_{t}S(t,x,h) = - \beta\big(\mu(t,x,h)\big) \displaystyle\frac{S(t,x,h)}{\mu(t,x,h)}\int_{O}K_\chi(x-y) I(t,y,h)\ud y \ ,\\[3mm]
\partial_{t}E(t,x,h) = \beta\big(\mu(t,x,h)\big)\displaystyle \frac{S(t,x,h)}{\mu(t,x,h)}\int_{O}K_\chi(x-y) I(t,y,h)\ud y - \theta E(t,x,h)\ ,\\[3mm]
\partial_{t}I(t,x,h) = \theta E(t,x,h) - \lambda I(t,x,h) - \delta I(t,x,h)\ ,\\[2mm]
\partial_{t}R(t,x,h) = \lambda I(t,x,h)\ ,\\[2mm]
\partial_{t}D(t,x,h) = \delta I(t,x,h)\ ,
\end{cases}
\end{equation}
For the motivation to our choice of the epidemiological model we refer to subsection \ref{sec:epimodel} and for a description of the economic model to subsection \ref{sec:MFG}. 

If the previous model is more handleable with respect to ones proposed in subsection \ref{subsec:dream} and subsection \ref{subsec:dream0}, it has the following structural problem. 
By \eqref{eq:sumM} we have
$$
\partial_t S(t,x,h)+\partial_tE(t,x,h)+\partial_t I(t,x,h)+\partial_tR(t,x,h)+\partial_tD(t,x,h)=0
$$
it follows 
$$
\partial_t\mu(t,x,h)=0.
$$
We conclude that the previous model is equivalent to choosing a stationary distribution $\mu(x,h)$. This remark lead us to study the model described in the following section.

\section{The hybrid MFG-SEIRD model}\label{sec:model}
In the following we introduce the SEIRD dynamic for the evolution of the epidemics, where the rate of transmission depends on the distribution of the population. Motivated by Section \ref{sec:listOfBadModels}, we assume that the distribution of people in space is given via a density function which is exogenous and does not change with time, but influences the incidence of the disease. Our model is clearly not suitable for describing any phase of the epidemic evolution, and this is even more true for the COVID-19 pandemics where measures as lockdown or total absence thereof significantly impact the time-evolution of the epidemic situation. Indeed we take population density in space as exogenous and given a priori; this means that we do not consider movement of people in time nor feedback effect of the epidemics dynamics on the population density. Therefore our model describes an epidemic phase in which the population density can be considered stationary in space-time; this corresponds for example to an initial stage of the spread of the disease, when regulators have not yet taken measures to prevent the propagation of the infection thus the population distribution is not affected by the dynamics of the epidemics, as well as a further stage in which movement-restricting measures have been lifted. In both cases the underlying idea that justifies considering some stationary distribution in space is that the mechanism that regulates people's movement has a time-scale that is significantly longer than the time-scale of the pandemic period we observe.

Thence, we have to face the issue of what population distribution to choose. There are of course many possibilities; the uniform distribution is surely the easiest choice but is not interesting in this regards as it would correspond to no differences across space, thus making the spatial extension of the model vanish in practice. Of course one could in principle test any distribution, but we are interested in using a distribution that has a solid economic justification.

Specifically, we choose as $\mu$ the marginal in the space variable of the stationary solution of a Mean-Field Game (MFG hereafter) describing equilibria in the space-time evolution of human capital subject to a utility maximization problem. Details are given in subsection \ref{sec:MFG}.
Considering  the space marginal of the solution to the FKP equation means in other word that we do not take into account the influence of the accumulation of human capital on the evolution of the epidemics. Up to our knowledge, the literature  has been mostly focusing more on the effect of the evolution of epidemics on the level of human capital \cite{Blea2010, DeShLiLi2021} instead of on the influence of such variable on the incidence term of the epidemiological model
A more realistic approach would consider non-stationary solutions of such MFG, leading to $\mu$ depending also on time; this is however not consistent with what argued in subsection \ref{subsec:dream00}. Another possibility would be to directly couple the mean-field game dynamics with our SEIRD dynamics, following the literature that discuss changes in the human capital distribution due to the epidemics dynamics; we refer to subsection \ref{subsec:dream} for some insight on the difficulties which might arise and we remark that studying this scenario is currently out of reach, both analytically and numerically.

\subsection{The SEIRD dynamics}\label{sec:epimodel}

We introduce now our epidemiological model. The space domain here is $\S^1$, the 1-dimensional torus. Positions in $\S^1$ are in one-to-one correspondence with positions in $[0,1)$; for every $x\in\bR$ the quantity $x\pmod{1}$ equals the fractional part of $x$, that belongs to $[0,1)$; the modulo $1$ operation allows us to consider the time evolution of positions as they belong to $\bR$ (which is more convenient for stochastic calculus tools) and then project them on the torus.\\
Moreover we need a consistent way to measure distances on the torus; we thus define, for every $x,y\in[0,1)$, the distance
\begin{equation*}
  d_{\S^1}(x,y):=\min\left\{(x-y)\pmod{1},(y-x)\pmod{1}\right\}\ .
\end{equation*}
This clearly corresponds to the arc-length distance on $\S^1$.

The population has constant size and is spatially distributed on $\S^1$ according to a density $\mu(x)$.

We assume that $\S^1$ has length $1$, and we mathematically represent it as the segment $[0,1]$ with the two endpoints identified when convenient for our purposes.

In our model we allow for a latency period between the time at which an individual gets infected and the time at which it can infect others; we do not consider births and deaths for ordinary causes, and we focus on the deaths caused by the infection. This implies that in the epidemic phase described by our model the number of deaths not ascribable to the infection is negligible with respect to the number of deaths in the population due to the infection. We also assume that all newborn individuals can contract the disease and that all recovered individuals are immune.

We divide the population into five non-intersecting classes, namely
\begin{itemize}
\item the class $S$ of \emph{susceptible} individuals, i.e. those that are currently not infected and can be infected;
\item the class $E$ of \emph{exposed} individuals, i.e. those that have got the infection but are not infectious yet, meaning thay cannot propagate;
\item the class $I$ of \emph{infectious} individuals, i.e. those that can spread the disease;
\item the class $R$ of \emph{recovered} individuals, i.e. those who healed and are thus immune to the disease;
  \item the class $D$ of \emph{dead} individuals, where death is due only to the disease and not to other causes.
  \end{itemize}
  Each variable $S, E,I, R, D$ is a function of time and space and represents the fraction of the total population that belongs to the corresponding class.

The main novelty we propose is the space-distribution dependent way in which we model the incidence of the disease. In the classical SIR-type models the form of the incidence is heuristically motivated as follows: if the population has size $N$, in a unit of time a single infectious individual is on average in contact with $C=cN$ individuals, for some proportionality constant $c>0$; among those individuals, $S/N$ are susceptible, thus leading to $cS$ contacts with susceptible individuals in one unit of time. Since transmission of the disease is not certain, only $pcS$ contacts produce new infectious, with $p\in(0,1)$. This must be multiplied by the number of infectious, leading to the incidence being $pcSI$, usually written as $\beta SI$. Since we take into account the distribution of individuals in space, this has to be modified accordingly: at the simplest level this implies substituting $N$ with $\mu (x)$, that is the number of people in the position $x$ of the infectious individual; this would lead again to the incidence having the form $\beta SI$. However there is evidence in the literature (see for example \cite{HuNiEc2013, LiSt2012, He2000}) that it is often not accurate to assume the average number of contacts $C$ of an infectious individual be linear in $N$; if we assume that it is not, then the above argument yields for the incidence term the form
\begin{equation}
  \label{eq:beta_1}
  \beta(\mu(x))\frac{S(t,x)}{\mu(x)}I(t,x)\ .
\end{equation}
Since the spatial structure of our model is continuous, with an incidence of the form (\ref{eq:beta_1}) one would have a contagion effect only with individuals that are \emph{exactly} at the same location at the same time. This is clearly not realistic, as the spread of the infection from an infectious individual occurs in some small area around her, with people being closer having higher probability of contracting the infection from the contagious individual, and people sufficiently far away having zero probability thereof. To model this feature we take the convolution in space of the variable $I$ with a smooth compactly supported symmetric kernel $K$; therefore the incidence at time $t$ and location $x$ depends not only on the number of infectious individuals at the same time and location, but also on the number of infectious individuals \emph{around} location $x$, with a maximal range of infection that is given by one half time the width of the support of $K$, and with infection strength that is maximal at location $x$ and decreases as one moves farther away from $x$. So, we propose the following form of the incidence:
\begin{equation}
  \label{eq:beta_ok}
  \beta(\mu(x)) \frac{S(t,x)}{\mu(x)}\int_{\mathbb{S}^1}K_\chi(x-y) I(y)\ud y
\end{equation}
where $\beta\colon\bR\to\bR$ is an increasing function such that $\beta(0)=0$ and $K_\chi(x)=\frac{1}{\chi} K\left(\frac{x}{\chi}\right)$ for $K\colon[0,1]\to\bR$ a smooth function, symmetric with respect to the point $\frac12$ and such that $K(x)=0$ for every $x\in\left[0,\epsilon\right]\cup\left[1-\epsilon,1\right]$, for some $\epsilon\in\left(0,\frac12\right)$. In this way $K$ can be identified with a function on $\S^1$ whose support is a strict compact subset of $\S^1$. For analogous approaches we refer to \cite{CoGaMaRo2020, TaHa2022}.

Denoting  by $\lambda$ the recovery rate of infectious individuals, by $\frac1\theta$ the latency period of the infection (i.e. the amount of time that passes between contracting the disease and becoming infectious, thus contributing to the spread of the disease), and by $\delta$ the death rate, our space-time dynamical model thus can be written as
\begin{equation}\label{eq:modelloDefinitivo}
\begin{cases} 
\partial_{t}S(t,x) &= - \beta\big(\mu(x)-D(t,x)\big) \displaystyle\frac{S(t,x)}{\mu(x)-D(t,x)}\int_{\Omega}K_\chi(x-y) I(y)\ud y \ ,\\[3mm]
\partial_{t}E(t,x) &= \beta\big(\mu(x)-D(t,x)\big)\displaystyle \frac{S(t,x)}{\mu(x)-D(t,x)}\int_{\Omega}K_\chi(x-y) I(y)\ud y - \theta E(t,x)\ ,\\[3mm]
\partial_{t}I(t,x) &= \theta E(t,x) - \lambda  I(t,x) - \delta I(t,x)\ ,\\[2mm]
\partial_{t}R(t,x) &= \lambda I(t,x)\ ,\\[2mm]
\partial_{t}D(t,x) &= \delta I(t,x)\ ,
\end{cases}
\end{equation}

together with suitable initial conditions and where, with respect to \eqref{eq:beta_ok},  in the incidence function we substracted from  $\mu(x)$ the  dead $D(t,x)$, since we expect deaths to have no influence in the incidence term.

For every $x\in\S^1$ we have
\begin{equation}\label{eq:sum}
  S(t,x)+I(t,x)+E(t,x)+R(t,x)+D(t,x) = \mu(x)\ ,
\end{equation}
while
\begin{equation*}
  \partial_{t}\left(S(t,x)+E(t,x)+I(t,x)+R(t,x)+D(t,x)\right)=0\ ,
\end{equation*}
confirming that the total population remains constant.
In this setting, contrary to some of the cited papers in the introduction, there is no movement of people due to the epidemic mechanism. 

In our numerical experience we  will consider  two main cases: $\beta$ in \eqref{eq:beta_ok} constant in $\mu(x)$, analyzed in \cite{PaLe2017} or $\beta$ having a nonlinear dependence in $\mu(x)$ as proposed in \cite{HuNiEc2013}. The case $\beta$  constant in $\mu(x)$ coincides with the  true mass action or proportioned mixing incidencewhere $\mu(x)$ replaces the total number of individual $N$. This case can be heuristically derived from a constant contact rate  not depending on the size of the population. Hence the per-link contact rate $\frac{C}{\mu(x)} $ decreases with larger density. However due to \eqref{eq:sum} an increase in $\mu(x)$ does not necessarily implies and increase in the number of susceptible individuals. We remark that, as argued in \cite{LiSt2012, He2000}, the case $\beta$ constant seem more consistent with the known result that daily contact patterns are  independent of community size.
On the other hand, in the present paper we focus also in the case where $\beta$ depends on $\mu(x)$ as in \eqref{eq:beta_ok} and in particular shows a nonlinear dependence as argued in \cite{HuNiEc2013}.

For simplicity we normalize here the total population to $1$. This implies that $\mu$ is a probability density, so that
\begin{equation*}
\int_{\mathbb{S}^1} \mu(x)\ud x = 1\ ,
\end{equation*}; if, in a population of size $N$ distributed in space according to $\mu(x)$, one wants to write the dynamics for the actual numbers $S^{(N)}=NS,E^{(N)}=NE,I^{(N)}=NI,R^{(N)}=NR,D^{(N)}=ND$ of individuals in each class one finds, for example in the first equation of (\ref{eq:modelloDefinitivo}),
\begin{equation*}
  \partial_{t}S^{(N)}(t,x)=- \beta^{(N)}\big(\mu^{(N)}(x)-D^{(N)}(t,x)\big) \displaystyle\frac{S^{(N)}(t,x)}{\mu^{(N)}(x)-D^{(N)}(t,x)}\int_{\mathbb{S}^1}K_\chi(x-y) I^{(N)}(y)\ud y,
\end{equation*}
where $\mu^{(N)}(x)=N\mu(x)$ and $\beta^{(N)}(z)=\beta\left(\frac{z}{N}\right)$. The other equations in (\ref{eq:modelloDefinitivo}) can be derived similarly.\\
The change of variables gets slightly more complicated if we allow for demography, i.e. rates of births and of deaths not due to the disease. We choose not to consider this extension here to simplify our exposition and presentation of numerical results. 

\subsection{The Mean Field Game model for the economy}
\label{sec:MFG}
We now introduce the MFG  whose stationary solution provides the exogenous spatial distribution for our epidemiological model. Such MFG arises as a continuum macroscopic description of a discrete (microscopic) multi-agent system in which agents interact through their human capital and decide how to move in space and their level of investment in professional skills in order to maximize a given utility function. An important feature of the economic model we consider is the presence of spatial spillovers on the accumulation of human capital and on consumption, where by spillover one denotes the influence of the presence of other individuals on the variable under consideration. The model  is therefore concerned with studying the effects on the abilities of each individual, caused by the proximity of other individuals. 

The MFG system is symmetric, in the sense that  both the interaction function and the utility function are the same for each agent, and agents are thus exchangeable. A key feature if that of  weak interaction, that is each agent has no influence on the other choices but the overall distribution has an impact on the optimization problem of each agents. This effect is encoded by a spatial interaction term consisting on a weighted average of the human capital of nearby agents, that is, a \textit{spatial spillover} on consumption. The MFG system consists in  two partial differential equations: an Hamilton-Jacobi-Bellman (HJB hereafter) equation describing the optimal control problem of each agent and a Fokker-Planck (FKP hereafter) or Kolmorogov equation describing the (optimal) evolution of the distribution of the population.

Heuristically, solutions to the mean-field game are approximate Nash equilibria for the model with a finite (large) number of agents. 
Well-posedness of the mean-field game and the numerics developed to find a stationary solution have been discussed in \cite{GRZ22}; here we describe the optimization problem for the system of interacting agents, we link it formally to the mean-field game and briefly comment on the numerical methods used to find $\mu$.


\subsubsection{The discrete multi-agent model}

In what follows we denote by $\sP_2$ the space of Borel probability measures on $\S^1\times\bR_+$ with finite second moment. Throughout this paper one can always assume that probability measures are absolutely continuous with respect to the Lebesgue measure on $[0,1]\times\bR_+$, implying that they have a density.

We consider a finite number $N$ of agents, indexed by $j=1,\dots, N$, whose positions and human capitals at time $t\geq 0$ are denoted respectively by $x_j(t)\pmod{1} \in \S^1$ and $h_j(t)\in\bR_+$, $j=1,\dots,N$.

Each agent has control over her own position, choosing the control process $v_j(t)$ in the equation
\begin{equation}
\label{eq:evolxnew}
\begin{cases}
  \ud x_j(t) = v_j(t) \ud t  +\epsilon \ud Z_j(t)\quad\text{ for }t\geq 0,\\
  x_j(0)=x_{j,0},
\end{cases}
\end{equation}
where
\begin{enumerate}[label=$(\roman{*})$]
\item $\epsilon>0$ is the noise intensity,
\item  $Z_j$ is a standard $1$-dimensional Brownian motion and
  \item $x_{j,0}$ are random variables taking values in $\S^1$, independent and identically distributed with some distribution $\mu_0^{1}$, that we assume to be absolutely continuous and to have finite second moment.
\end{enumerate}
The spatial dynamics of the agents affect their human capital, which evolves according to the equation
\begin{equation}
  \label{eq:evolhnew1}
  \begin{cases}
    \ud h_j(t) = f_j(t) h_j(t)^\alpha \bar h_{-j} (t)^{\xi} -\zeta h_j(t)\ud t +\sigma h_j (t) \ud W_j(t)\quad\text{ for }t\geq 0,\\
    h_j(0)=h_{j,0},
  \end{cases}
\end{equation}
where
\begin{enumerate}[label=$(\roman{*})$]
\item $f_j(t)\in[0,1]$ is the ``investment in education'' and can be chosen by the agent;
\item $ \alpha,\xi \in (0,1)$ are constants such that $\alpha>\frac12$ and $\alpha + \xi <1$; the motivation for introducing $f$ and $\alpha$ is that the human capital is partly utilized to increase the human capital itself ($f_j(t)h_j(t)^\alpha$ where $h_j(t)^\alpha$ is the ``personal income'') and partly to produce final goods ($(1-f(t))h_j(t)$, see (\ref{eq:funzionale1new}) below); 
\item $\zeta>0$ is the depreciation rate of human capital;
\item $\sigma>0$ is the intensity of the noise, that is here chosen as a multiplicative noise to ensure positivity of the solutions;
\item $W_j$ is a standard $1$-dimensional Brownian motion;
  \item $h_{j,0}$ are random variables taking values in $\bR_+$, independent and identically distributed with some distribution $\mu_0^{2}$, that we assume to be absolutely continuous and to have finite second moment;
\item $\bar h_{-j} (t)$ is an average depending on $x_j(t)$:
  \begin{equation}
    \label{eq:barhbisnew}
    \bar h_{-j} (t) := \frac{\sum_{k\neq j} \eta \left(d_{\S^1}\left(x_k(t) - x_j(t)\right)\right) h_k(t)}{\sum_{k\neq j} \eta \left(d_{\S^1}\left(x_k(t) - x_j(t)\right)\right)},
  \end{equation}
where the function $\eta$ encode the interaction between agents and is given by
\begin{gather*}
  \eta\colon\left[0,1/2\right]\to\bR_+\\
  \eta(x)=1 \mbox{ if } x\leq \epsilon_1,\quad   \eta(x) =\epsilon_2 \mbox{ if } \frac{1}{2}-\epsilon_3\leq x\leq \frac{1}{2},
\end{gather*}
and $\eta(x)$ is assumed linear between $\eps_1$ and $\frac{1}{2}-\eps_3$ so that the function results Lipschitz continuous in all its domain $[0,1/2]$,
where $\epsilon_1<\frac{1}{2}, 0<\epsilon_2\ll1, \frac{1}{2}-\epsilon_3 >\epsilon_1$.
\end{enumerate}

All initial conditions $x_{j,0}$ and $h_{j,0}$ and the Brownian motions $Z_j$ and $W_j$ are defined on some common probability space $(\Omega, \sF,\bP)$ and are supposed to be independent.

The function $\eta$ is simply a smooth version of the indicator function of a right-neighborhood of $0$; indeed $\eta \left(d_{\S^1}\left(x - y\right)\right)$ equals $1$ whenever $x$ and $y$ are closer than $\eps_1$, and equals $0$ if $x$ and $y$ are far from each other. Therefore the interaction term $\bar h_j$ should be understood as follows: at each time $t$, agent $j$'s human capital benefits mostly from that of all other agents that at time $t$ are located within $\eps_1$ from her. The value of human capital coming from this interaction term is divided by the number of agents that are located within $\epsilon_1$ from her; this mean-field-type scaling is necessary to ensure that the interaction term does not blow up when the number of agents increases.




Both the equations for $x_j$ and $h_j$ depend on controls (respectively $v_j(\cdot)$ and $s_j(\cdot)$); for each $j$, the aim of agent $j$ is to choose such controls as to maximize the gain functional
\begin{equation}
  \label{eq:funzionale1new}
  J_{j}\left(x_{0,j},h_{0,j};v_j(\cdot), f_j(\cdot)\right):=\mathbb{E} \left [ \int_{0}^{+\infty} e^{-\rho t}\left(u\left( [(1-f(t)) h_j(t)^\alpha]^{1-\gamma} \bar h_{-j}(t)^\gamma A\left(x_j(t)\right)\right)  -a(v_j(t))\right) dt  \right ],
\end{equation}
where $\bE$ denotes expectation with respect to $\bP$. The parameters and functions appearing in the definition  of $J$ are:
\begin{enumerate}[label=$(\roman{*})$]
\item $\rho>0$ is the discount factor;
  \item $u$ is the utility function, which takes the constant relative risk adversion form
	\begin{equation}\label{eqn:utility}
	u(z)=\frac{z^{1-p}}{1-p}
      \end{equation}
      for some $p\in (0,1)$.
    \item $-a(v(t))$ represents the cost for the energy employed in the displacement;
    \item $[(1-f(t)) h_j(t)^\alpha]^{1-\gamma}$ represents the final goods produced thanks to the human capital;
    \item $ \bar h_j(t)^\gamma$ is the spatial spillover on consumption;
      \item $A(x)$ represents the local amenities; here we choose $A(x)=\frac{1}{2}\sin\left(2\pi\left(x-\frac{1}{4}\right)\right)+1$.
    \end{enumerate}
        

Each agent's choice of the controls determines the general dynamics of all agents, because of the interaction term $\bar h_{-j}$. 

\subsubsection{The Mean-Field Game system}
The multi-agent system described above is symmetric, in the sense that agents are indistinguishable in their contribution on each other's human capital and they are all maximizing the same game functional; thus it is reasonable to look for Nash equilibria of the system. Nash equilibria are however impossible to compute exactly for a system with a large number of agents. Therefore a standard way to investigate equilibria is to take the limit as the number of agents goes to infinity and to study the limit system. A solution to the limit system will then approximate a Nash equilibrium for the original multi-agent system when the number $N$ of agents is sufficiently large. The limit system is usually referred to as a mean-field game (MFG hereafter); although a mathematically rigorous proof of the connection between the multi-agent system and the MFG has been carried out only in simple cases, it now widely accepted in economics and mathematical modeling the idea that MFGs are the correct tool to study situations as the one we are describing here, see for example \cite{ABLLM14}.

To write down the MFG it is convenient to introduce the function $F\colon\S^1\times\sP_2\to\bR_+$ defined as
\begin{equation}\label{eq:barh}
  F(x, \xi)=\frac{\int_{\S^1\times\bR_+}\eta\left(d_{\S^1}(x,y)\right)k\xi\left(\ud y,\ud k\right)}{\int_{\S^1\times\bR_+}\eta\left(d_{\S^1}(x,y)\right)\xi\left(\ud y,\ud k\right)};
\end{equation}
in this way equation the interaction term $\bar h_{-j}(t)$ takes the form
\begin{equation*}
  \bar h_{-j}(t)=F\left(x_j(t),S(t)\right)
\end{equation*}
where $S(t)$ is the \emph{empirical measure} of the system, defined as the sum of Dirac deltas
\begin{equation*}
  S(t)=\frac1N\sum_{j=1}^N\delta_{\left(x_j(t),h_j(t)\right)}\ .
\end{equation*}

The MFG system associated to the optimal control problem \eqref{eq:funzionale1new} is the following system of a HJB equation coupled with a FKP equation, 
\begin{equation}\label{MFG}
\begin{cases}
 \rho V=H_1\big(x,h,\mu(t), \partial_hV\big)+ \frac12 \chi^2h^2 \partial^2_{hh}V+\frac12 \eps^2 \partial^2_{xx}V+H_0(\partial_xV)\ ,\\[3mm]
 \partial_t\mu=\frac{1}{2}\chi^2\partial^2_{hh}(h^2\mu)+\frac{1}{2}\eps^2\partial^2_{xx}\mu-\partial_x\left(\partial_p H_0(\partial_x V)\mu\right)-\partial_h\left(\partial_p  H_1(x,h,\mu(t), \partial_hV)\mu\right)\ ,
\end{cases}
\end{equation}
whose solution is a couple $(V,\mu)$ of functions of the three variables $t>0$, $x\in\S^1$ and $h>0$ such that for every $t>0$ $\mu(t,\cdot,\cdot)$ is a probability density on $\S^1\times\bR_+$. The function $V$ solves the first PDE above in the classical sense, while the function $\mu$ solves the second PDE above in the sense of distributions. The system is equipped with an initial condition for $\mu$: we require that
\begin{equation*}
  \mu(0,x,h)=\mu_0(x,h)=\mu_0^{(1)}(x)\mu_0^{(2)}(h) \text{ for every }(x,h)\in\S^1\times\bR_+\ ,
\end{equation*}
Moreover we need a set of boundary conditions for the $h$-boundary $\{h=0\}$; we impose Dirichlet boundary conditions, namely
\begin{equation*}
   V(t,x,0)= 0,\quad \mu(t,x,0)=\bar \mu(t,x) \text{ for every }(t,x)\in\bR_+\times\S^1
\end{equation*}
where for every $t>0$ $\bar \mu(t,\cdot)$ is a bounded probability density on $\S^1$ with finite second moment; this forces the compatibility condition
\begin{equation*}
 \mu_0(x,0)=\bar \mu(0,x).
\end{equation*}
The function s $H_0$ and $H_1$ appearing in (\ref{MFG}) are the Hamiltonians of the two coupled control problems in the multi-agent system, one for the control $v$ and the other for the control $s$. $H_0$ is defined as
\begin{gather*}
  H_0\colon \bR\to\bR\\
  H_{0}(p)=\sup_{v \in K}\left \{-a(v)+vp\right\};
\end{gather*}
$H_1$ instead is defined as
\begin{gather*}
  H_1\colon\bR_+\times\S^1\times\bR_+\times\sP_2\times\bR\to\bR\\
H_1(t, x,h,\xi, q)=\sup_{ f\in [0,1]}\left\{(fh^\alpha F(x,\xi)^\xi-\zeta h) q+u\left( A(x)[(1-f) h^\alpha]^{1-\gamma} F(x,\xi)^\gamma\right)\right\},
\end{gather*}
For every fixed function $\bR_+\ni t\mapsto \mu(t)\in\sP_2$ the first equation in (\ref{MFG}) is the HJB equation for the optimization problem of a representative agent that interacts with a continuum of agents whose locations and human capitals are distributed according to $\mu$; the value function for such optimization problem is then a solution of such HJB equation. The second equation in (\ref{MFG}) instead describes the time evolution of the space and human capital distribution of a continuum of agents, given a solution $V$ of the first equation. The MFG consists precisely of the coupling of the two.

As explained in the introduction of the paper and in the introductory paragraph of Section \ref{sec:model}, we use as our a priori spatial distribution the $x$-marginal of a stationary solution $\mu$ to the second equation of (\ref{MFG}). Such marginal is a probability density on $\S^1$.





\subsubsection{Numerical solution to the steady state problem}\label{sec:numericalMFG}
\begin{figure}[t!]
\includegraphics[width=\textwidth]{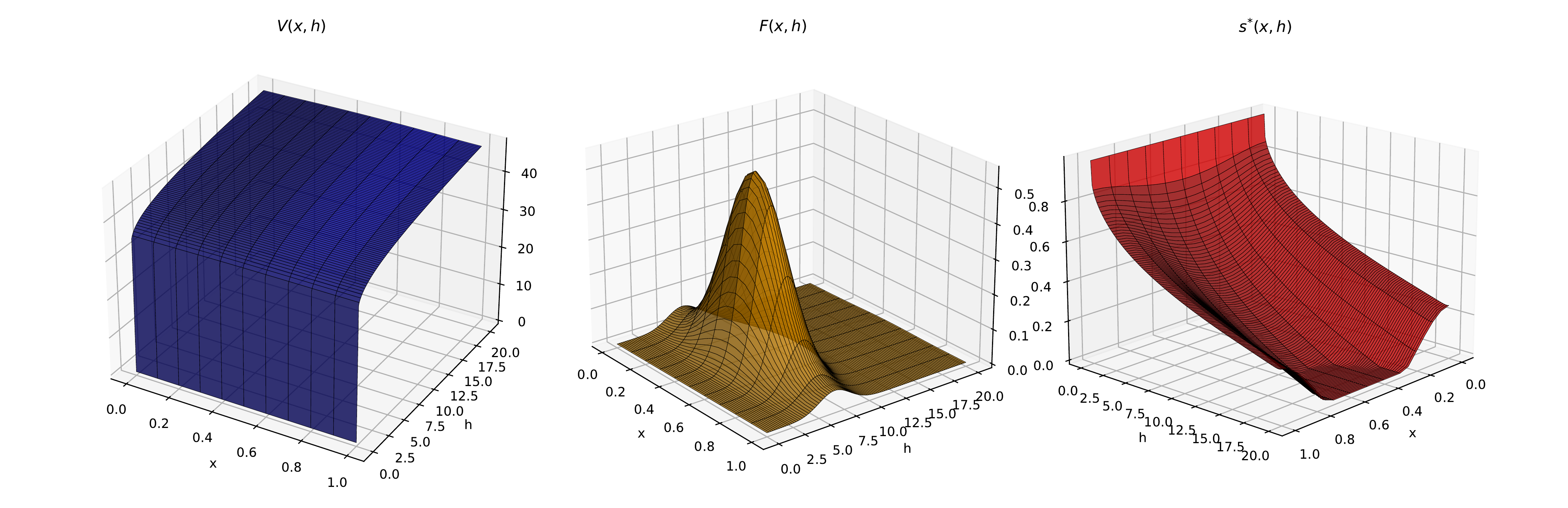}
\caption{Numerical equilibrium solution to system of equation \eqref{eq:MFGsteadyState}. The left panel shows the plot of the value function $V(x,h)$ (blue), the one in the middle that of the Fokker-Plank (right) and right most one shows the optimal investment (red). We can see in the central panel how the bulk of the population concentrates more in those areas corresponding to higher spatial amenities ($x \approx 0.5$), as well as around an intermediate value ($h \approx 5$) along the $h$ variable. The optimal investment rate exhibit a decreasing behaviour with respect to the $h$ variable, and even more reaches the value of zero for in the area of the domain corresponding to highly educated individuals living within the central large city.}
\label{fig:MFG}
\end{figure}

In this section we study the steady state of the system of equations \eqref{MFG} from a numerical view point. Therefore we consider the following problem
\begin{equation}\label{eq:MFGsteadyState}
\begin{cases}
\rho V=H_1(x,h,\mu, D_hV)+ \frac12 \chi^2h^2 D^2_{hh}V+\frac12 \eps^2 D^2_{xx}V+H_0(D_xV) &\mbox{ in } \mathbb{S}^1\times \R_{++}\\
\frac{1}{2}\chi^2\partial^2_{hh}(h^2\mu)+\frac{1}{2}\eps^2\partial^2_{xx}\mu = \partial_x\left(D_p H_0(D_x V)\mu\right)+ \partial_h\left(D_p  H_1(x,h,\mu, D_hV)\mu\right)& \mbox{ in }  \mathbb{S}^1\times \R_{++}
\end{cases}
\end{equation}
with Dirichlet boundary conditions as in the non stationary case 
\begin{equation*}
V(x, 0)= 0,\quad \mu(x,0)=\bar \mu(x) \quad  \mbox{ in } \mathbb{S}^1.
\end{equation*}
The problem of finding numerical solutions to MFGs is in general very challenging. For a full discussion of the methodology applied to find numerical solutions to \eqref{eq:MFGsteadyState} we refer to \cite{GRZ22}. In this manuscript we limit ourselves to a brief outline of the main elements. 
The numerical discretization is performed by means of classical finite difference method. However, for MFGs problem one has to pay special attention on the way the discretization of the Hamiltonian is performed. In particular in order to ensure consistency of the discretization scheme we adopted the \emph{Kushner-Dupuis} approach, see \cite{achdou2010mean} and also \cite[Chapter 2.]{achdou2020mean} for a more general summary on the topic. 
The choice of an appropriate truncation of the domain $\RR_{+}$ for the $h$ variable is another difficult topic that one has to take care of. Numerical simulations in fact are performed on the truncated domain $S^1 \times [0,H_{max}]$ where the upper  bound value $H_{max}$ has to be imposed artificially. In particular, we had to adopt an iterative strategy in the choice of $H_{max}$ by computing a sequence of solutions with increasing $H_{max}$ up until the bulk of the distribution of $\mu$ was fully enclosed in the truncated domain. See Table \ref{tab:paramMFG} for a full list of the parameters used, including the threshold $H_{max}$. 
\begin{table}[]
  \begin{tabular}{rlrll}
  \textbf{Parameter}        & \multicolumn{1}{l|}{\textbf{Value}}                               & \textbf{Parameter}   & \textbf{Value} & \textbf{} \\ \cline{1-4}
  $L$ (length of the torus) & \multicolumn{1}{l|}{$1$}                                          & $\alpha$             & $0.5$          &           \\
  $R$ (radius for $\eta$)                       & \multicolumn{1}{l|}{$0.3$}                                        & $\xi$                & $0.1$          &           \\
  $H_{max}$                 & \multicolumn{1}{l|}{$15$}                                         & $\gamma$             & $0.15$         &           \\
  $A(x)$                    & \multicolumn{1}{l|}{$\frac{1}{2}\sin (2\pi (x-\frac{1}{4})) + 1$} & $\chi$             & $0.7$          &           \\
  $\epsilon$                & \multicolumn{1}{l|}{$0.5$}                                          & $\zeta$             & $0.15$         &           \\
  $p$                    & $0.1$                                                             & \multicolumn{1}{l}{} &                &          
  \end{tabular}
  \caption{Parameters for system of equation \eqref{eq:MFGsteadyState} Figure \ref{fig:MFG}.}
  \label{tab:paramMFG}
  \end{table}
One of the key parameters characterizing the stationary solution of the MFG is the function $A(x)$, $x \in S^1$, describing the distribution over space of local amenities. 
In fact, in the stationary case, the only source of spatial heterogeneity is encoded in the exogenous function $A(x)$. This is not true in general when studying the time evolution of the system out of equilibrium, in which the initial distribution plays an important role into determining the spatial distribution of agents at any given period of time $t$. However, when studying the system at equilibrium the influence of the initial configuration vanishes as $t$ grows to $+ \infty$ and the only remaining factor affecting the shape of the solutions has to be accounted to the function $A$. In particular, one can see that in absence of heterogeneity, i.e. if $A(x)$ is uniform, solutions to system of equations \eqref{eq:MFGsteadyState} are translation invariant on the space variable. That is if one assumes that $\left(V(x,h),\mu(x,h)\right)$ is a solution to the stationary problem, the same holds for $\left(V(x+x',h),\mu(x+x',h)\right)$ for any $x' \in S^1$. On the contrary, when $A(x)$  is not uniform, the spatial distribution of agents is subject to a positive influx from local exogenous amenities. In particular one can imagine that, being the objective functional \eqref{eq:funzionale1new} monotone increasing with respect to the function $A$, each agent will relocate trying to obtain the best possible level (net of reallocation costs) of local amenities, that is where $A(x)$ is higher. Figure \ref{fig:MFG} shows the equilibrium solution to \eqref{eq:MFGsteadyState} when the function $A$ has been taken to have a single maximum peak at $x = 0.5$. We see in particular how agents distribute following the same shape as argued in the previous lines. See moreover Table \ref{tab:paramMFG} for a full list of parameters, including the function $A$.

The case shown in Figure \ref{fig:MFG} is precisely the case we had in mind, and that we will also employ in the next section, while discussing the spatial epidemic model \eqref{eq:modelloDefinitivo}. In particular, we imagine that having $A(x)$  with one single peak, represents the presence of a single and large city whose center is located in correspondence of the maximum of $A$. This idea is reinforced by the shape of $\mu(x,h)$ as shown in Figure \ref{fig:MFG} where the bulk of the population concentrates in the middle. 

\section{Numerical results on the spatial SEIRD model with heterogeneous distribution}\label{sec:numericaSEIRD}
In this section we present some qualitative results related to the SEIRD model \eqref{eq:modelloDefinitivo} when the spatial distribution of agents is taken as the space marginal of the stationary distribution \eqref{eq:MFGsteadyState} . We will see, with the aid of numerical simulations, how the role of the population density plays a crucial role in the evolving of the epidemic through time. 

The choice of parameters plays a key role in the evolution of epidemiological models. From  the theoretical point of view different choices may lead do different qualitative behaviour in the equilibrium, exhibiting a phenomenon of phase transition for the steady state. At the same time, different sets of parameters describe also specific characteristics of different diseases, as well as of that of the same disease in different countries.  Several works, many of which appeared  in the recent years, focus on the problem of statistical estimating those parameters from available datas. 
\begin{figure}
\centering\includegraphics[width=5.5cm]{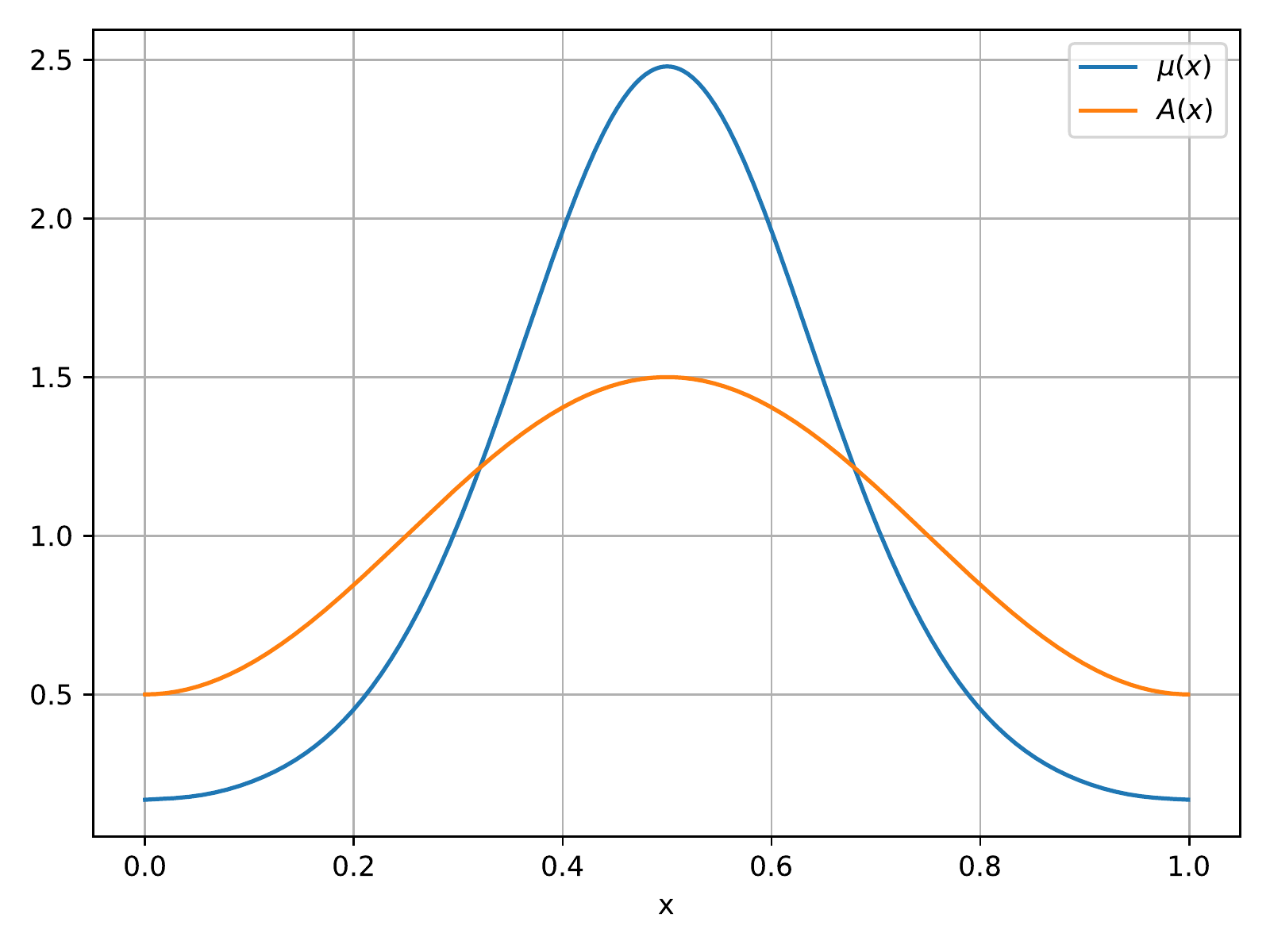}
\caption{The function $A(x)$ (orange) vs. the spatial distribution of agents $\mu(x)$ (blue) obtained from the Mean Field Game. We see how the spatial distribution of agents resembles the shape of $A(x)$. This accounts to the fact that the objective function \eqref{eq:funzionale1new} is monotone increasing in the function $A$ with respect to the position of agents, and reflect the tendency of agents to concentrate more where local amenities are more present.}\label{fig:muA}
\end{figure}


In the present paper we don't aim to represent a real-world scenario, since we are still limited by the current level of mathematics and numerical capabilities, but to bring a methodological contribution. However, in order to be as close as possible to reality, we decided to consider the case of COVID-19 for the US, hence considering a set of epidemiological parameter as close as possible to the scenario we have in mind. As mentioned above, a plethora of papers are present on the topic of parameter estimation, but we take the work \cite{Korolev2021} as a reference, therefore taking $(\beta, a,\gamma,\delta) = (0.9,0.25,0.075,0.025)$ (from the case $\sigma = 1/4, \gamma = 1/10$ in \cite{Korolev2021}). See moreover Table \ref{tab:paramEpidemic} for a full list of parameters.
\begin{figure}
\includegraphics[width=\textwidth]{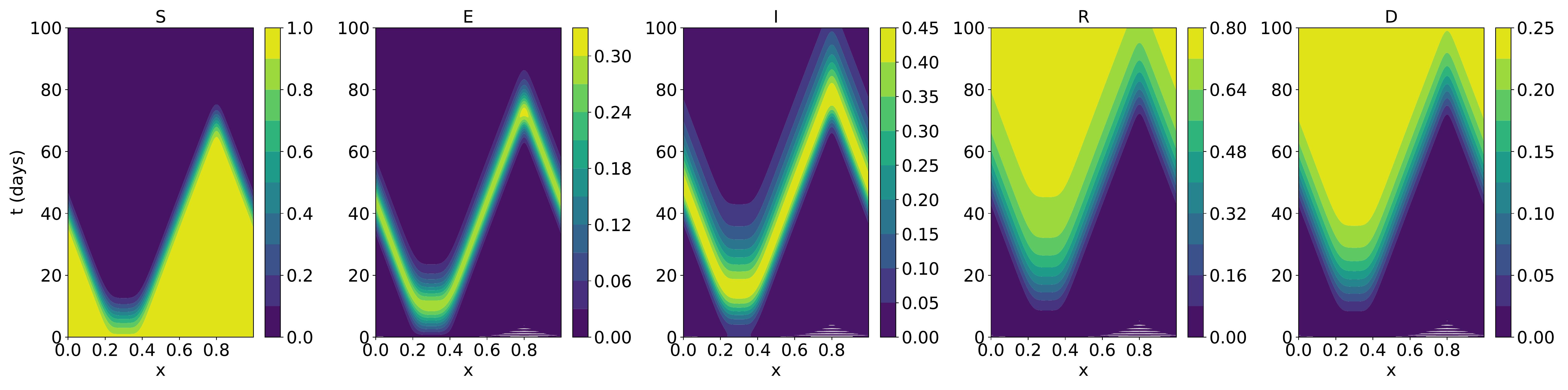}
\caption{Space-time evolution of the SEIRD model \eqref{eq:modelloDefinitivo} with uniform population distribution $\mu(x) \equiv 1$. The horizontal axis represents the space variable, while time is shown vertically. We see that the infection originates around $x = 0.3$ and briefly spreads to neighboring locations. The farthest locations from where the infection outbreaks are the last one to receive the contagious wave of infected, around $x \approx 0.8$. The final configuration approaches the uniformity in space in terms deaths. The value of $\beta$ is kept homogeneous as in Table \ref{tab:paramEpidemic}}
\label{fig:constantMUconstantBeta}
\end{figure}
As previously mentioned in subsection \ref{sec:MFG} the shape of $A$ dictates the population distribution along the space variable. In Figure \ref{fig:muA} we clearly see this fact, by noticing the shape of the function $\mu(x)$. 
Since the bulk of the mass of population concentrates around the point $0.5$, in order to better appreciate the effects of the population distribution, we take an initial configuration for the SEIRD model \eqref{eq:modelloDefinitivo} where the outbreak of the epidemic origins from the outskirts. Hence, we set the initial condition to have a small cluster of infects centered around the point $0.3$ in the class $I$. 
Namely we set
\begin{equation*}
S(0,x) = \mu(x)-I_{0}\mathds{1}_{B^{\eps}(0.3,r_{0,I})},\quad E(0,x) = 0,
\end{equation*}
\begin{equation*}
I(0,x) = I_{0}\mathds{1}_{B^{\eps}(0.3,r_{0,I})},\quad R(0,x) = 0,\quad D(0,x) = 0,\quad \forall x \in S^1
\end{equation*}
where $\mathds{1}_{B^\eps(0,r)}$ is a suitable mollification of the indicator function in order to have a smooth initial condition, and the parameters $I_{0}$ and $r_{I,0}$ are set to the values $0.01$ and $0.1$ respectively. Notice that the parameter $\epsilon > 0$ is just used as a regularization parameter for the indicator function. Also, it can be taken arbitrary small since the function $\mathds{1}_{B^\eps(0,r)}(x)$ is of class $C^\infty$ as soon as the parameter $\epsilon$ is positive, without any further assumption on it's magnitude. Of course many different initial condition are possible. The one that we selected is chosen so that the initial cluster, center in $x = 0.3$, is not located in the location of maximum density, $x = 0.5$. Moreover, being the model defined on the $1$-dimensional torus, taking initial cluster of infected to the farthest point from the peak of the agent distribution would make the euclidean distance between the two points identical moving in either of the two possible directions (left or right). Therefore, we selected an initial condition which tries to avoid any kind of misleading results due to spatial symmetry.

Concerning the interaction kernel $K_{\chi}$ in equation \eqref{eq:modelloDefinitivo} we recall that it describes the fact that people may get in contact not only with those individual located in their same position, but also from the neighboring areas. In particular we consider the function $K$ as a suitable mollification\footnote{The mollification has no modellistic meaning and it is only here to ensure regularity of the solution without any technical difficulty.} of the \emph{hat function} $H(x)$, defined as
\begin{equation*}
H(x) = \begin{cases}
1 - \abs{x} &\text{ if } \abs{x} < 1,\\
0 &\text{otherwise.}
\end{cases}
\end{equation*}
Then we adopt the following rescaling 
\begin{equation*}
K_{\chi}(x) = \frac{1}{\chi}K\left(\frac{x}{\chi}\right)
\end{equation*}
that keeps the total area of the kernel untouched, while changing the radius at which infections can be transmitted. This way of defining the rescaling is fundamental in order to disentangle the effects of the parameter $\chi$, which measure the radius of interaction, from the function $\beta$ (parameter in the case it is taken constant), which measure the transmission rate. 
In all the cases that we consider we set $\chi = 0.04$. Since our model is defined on the periodic one dimensional space $S^1$ this value has no particular meaning in absolute terms. However, one important point to remark is that the selected value is chosen one order of magnitude smaller than the radius of interaction $R$ of the Mean Field Game system (see Table \ref{tab:paramMFG}). This reflects the fact that information, being in general an immaterial good, can travel much further (per unit of time) than viruses or people. Hence the exchange of information that agents receive, in the form of spatial spillovers from other individuals, can be seen at much higher distance that the radius of infection, therefore supporting a higher value of $R$ with respect to $\chi$. 
\begin{table}[]
\centering
\begin{tabular}{rl|rll}
\textbf{Parameter}      & \textbf{Value} & \textbf{Parameter}   & \textbf{Value} & \textbf{} \\ \cline{1-4}
$\theta$                     & $0.25$         & $\chi$                  & $0.04$         &           \\
$\lambda$                & $0.075$        & $I_0$                & $0.01$         &           \\
$\delta$                & $0.025$        & $r_{I,0}$            & $0.1$          &           \\
$\beta$ (when constant) & $0.9$          & $\beta(\mu)$         & $0.9\frac{\sqrt{\mu(1+\mu)}}{\sqrt{2}}$         &           \\
                        &                &                      &                &           \\
                        &                & \multicolumn{1}{l}{} &                &          
\end{tabular}
\caption{Parameters for system of equation}
\label{tab:paramEpidemic}
\end{table}
\begin{figure}
\includegraphics[width=\textwidth]{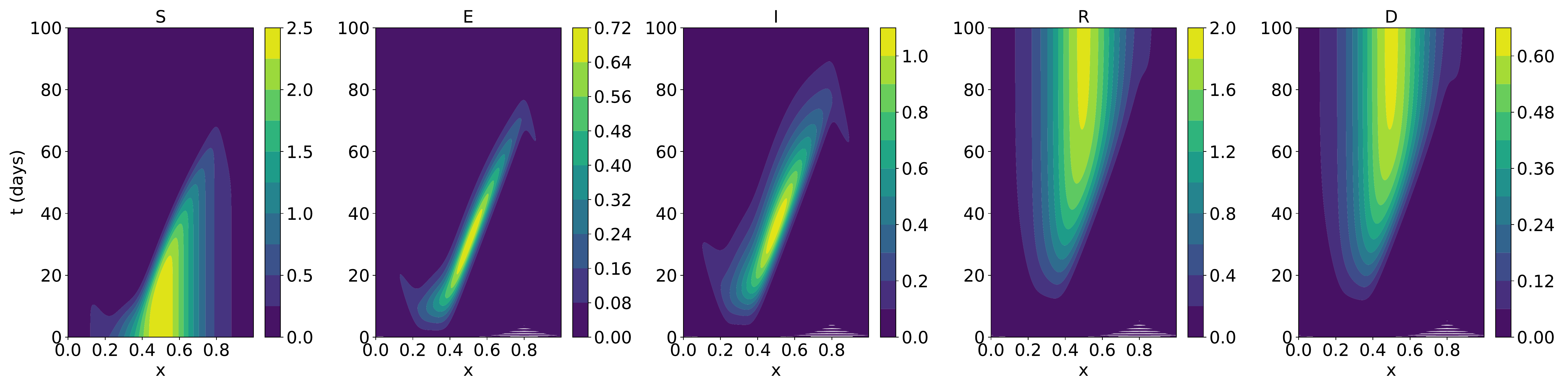}
\caption{Space-time evolution of the SEIRD model \eqref{eq:modelloDefinitivo} with spatial distribution $\mu(x)$ obtained from the Mean Field Game system \eqref{eq:MFGsteadyState}. The horizontal axis represents the space variable, while time is shown vertically. The infectious disease originates around $x = 0.3$. We can appreciate from the graph how the diffusion of the virus favours the locations where the total population is higher. The final configuration shows the same behaviour. The value of $\beta$ is kept constant as in Table \ref{tab:paramEpidemic} to the value $0.9$.}
\label{fig:MFGmuConstantBeta}
\end{figure}

\begin{figure}
\includegraphics[width=\textwidth]{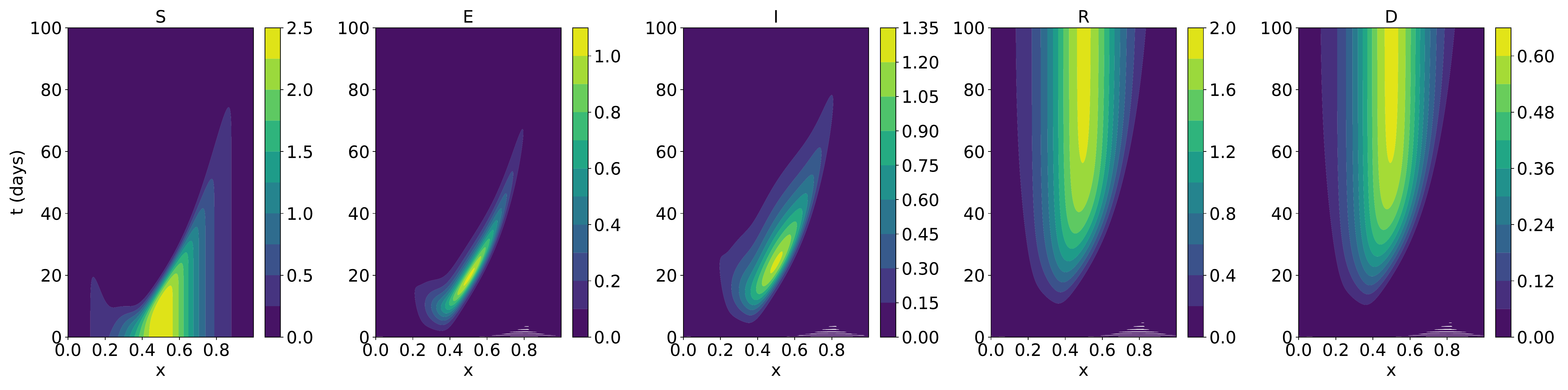}
\caption{Space-time evolution of the SEIRD model \eqref{eq:modelloDefinitivo} with spatial distribution $\mu(x)$ obtained from the Mean Field Game system \eqref{eq:MFGsteadyState}. The horizontal axis represents the space variable, while time is shown vertically. The infection rate $\beta$ is taken to be dependent on the population density, see Table \ref{tab:paramEpidemic}. With respect to the scenario with constant $\beta$ here we see the differences in the evolution infected. In particular the spreading of the disease is much more localized around the center of the graph, were not only the population density is higher, but also the infection rate increases. The opposing effect happens closer to $x = 0$ (analogously $x = 1$ due to periodicity). See in particular panels corresponding to classes $E$ and $I$.}

\label{fig:MFGmuSQRTBeta}
\end{figure}

We simulate the SEIRD model \eqref{eq:modelloDefinitivo} by means of classical finite difference method. We run the model for a time period of $100$ days, and analize the result in a space-time plot, see Figures \ref{fig:constantMUconstantBeta} and \ref{fig:MFGmuConstantBeta}. We first start by looking at the effect that the population density has on the evolution of the pandemic. Figure \ref{fig:constantMUconstantBeta} shows the case where the spatial distribution of population is neglected, and a completely flat distribution is taken, while the subsequent figure, Figure \ref{fig:MFGmuConstantBeta}, describes the same scenario but with population density as obtained from the Mean Field Game. In both cases $\beta$ is constant.
Finally we also consider the case where the infection rate $\beta$ is itself dependent on the density $\mu$, see Figure \ref{fig:MFGmuSQRTBeta}. This choice is not new, and has been explored from example in \cite{HuNiEc2013} taking into account the average contact rate the each individual has. The rationale behind it is the following: the number of contact has an average person has, has a direct impact of the characteristic infectious rate of a given disease (like COVID-19 as in our case). Therefore, it is reasonable to assume that, since the number of contact per individual increases on average, as a function of the population density, so does the infection rate $\beta$. Therefore, we introduce the function $\beta(\mu)$. Following \cite[p. 129]{HuNiEc2013} we take
\begin{equation*}
\beta(\mu) = 0.9\frac{\sqrt{\mu(1+\mu)}}{\sqrt{2}}.
\end{equation*}
The choice of the constants appearing in the function $\beta(\mu)$ has the following reasoning: the factor $\sqrt{2}$ is introduced so that, in the case of a uniform distribution for $\mu$ (in our case $\mu(x) \equiv 1$), then the right part of the formula reduces to the factor $1$. By this choice, the factor $0.9$ in front makes it so that when $\mu(x) \equiv 1$ the value of $\beta(\mu)$ coincides with the one used in the case where the distribution of agent is assumed uniform across the space, that is $\beta = 0.9$, see Figure \ref{fig:constantMUconstantBeta} for numerical simulations in that case. 

\section*{Declarations}
The authors have no competing interests to declare that are relevant to the content of this article.

\bibliography{biblio} 


\end{document}